\documentclass[12pt]{amsart}
\usepackage{amssymb}
\theoremstyle{plain}
\newtheorem{theorem}{Theorem}
\newtheorem{corollary}{Corollary}

\newtheorem{proposition}{Proposition}
\theoremstyle{example}
\newtheorem{example}{Example}
\theoremstyle{definition}

\theoremstyle{remark}

\numberwithin{equation}{section}
\def\K#1#2{\displaystyle{\mathop K\limits_{#1}^{#2}}}

\setbox0=\hbox{$+$}
\newdimen\plusheight
\plusheight=\ht0
\def\+{\;\lower\plusheight\hbox{$+$}\;}

\setbox0=\hbox{$-$}
\newdimen\minusheight
\minusheight=\ht0
\def\-{\;\lower\minusheight\hbox{$-$}\;}

\setbox0=\hbox{$\cdots$}
\newdimen\cdotsheight
\cdotsheight=\plusheight
\def\cds{\lower\cdotsheight\hbox{$\cdots$}}
\begin{document}
\title[Real Numbers With Polynomial Continued Fraction Expansions]
 {Real Numbers With Polynomial Continued Fraction  Expansions}
\author{J. Mc Laughlin}
\address{ Department of Mathematics\\
 Trinity College\\
300 Summit Street, Hartford, CT 06106-3100}
\email{james.mclaughlin@trincoll.edu}
\author{ Nancy J. Wyshinski}
\address{Department of Mathematics\\
       Trinity College\\
        300 Summit Street, Hartford, CT 06106-3100}
\email{nancy.wyshinski@trincoll.edu}

\date{January 21, 2004}
\begin{abstract}
In this paper we show how to apply various techniques and theorems
(including  Pincherle's theorem, an extension of Euler's formula equating
infinite series and continued fractions, an extension of the
corresponding transformation that equates infinite products and
continued fractions, extensions and contractions of continued
fractions and  the Bauer-Muir transformation) to derive
infinite families of in-equivalent polynomial continued fractions
in which each  continued fraction has the same limit.

This allows
us, for example, to construct infinite families of polynomial
continued fractions for famous constants like $\pi$ and $e$,
$\zeta{(k)}$ (for each positive integer $k\geq 2$), various special
functions evaluated at integral arguments and various algebraic
numbers.

We also pose several questions about the nature of the set of real
numbers which have a polynomial continued fraction expansion.
\end{abstract}
\maketitle

\section{Introduction}

A polynomial continued fraction is a continued fraction
  $ K_{n=1}^{\infty}\displaystyle{a_{n}/b_{n}}$ whose partial numerators and denominators
 are integers and where, for some tail of the continued fraction,
$a_{n}=f(n)$, $b_{n}=g(n)$, for some
polynomials $f(x),\,\,g(x) \in \mathbb{Z}[\,x\,]$.

One reason polynomial continued fractions are
of interest is that sometimes a famous number whose regular continued fraction expansion
does not have a predictable pattern can  be expressed as polynomial
continued fractions (so that the
partial quotients do have a predictable pattern).
For example the first continued fraction
giving values for $\pi$\,(due to Lord Brouncker, first published in
 ~\cite{W65})  is of this type:
\vspace{3pt}
\begin{equation}\label{E:lb}
\frac{4}{\pi} = 1 +
\K{n=1}{\infty}\displaystyle{\frac{(2n-1)^{2}}{2}}.
\end{equation}

Another reason polynomial continued fractions are interesting is that the partial quotients are integers
and standard theorems can sometimes be used to decide if a number is irrational. For example,
there is the following criterion given by Tietze,
extending the famous Theorem of Legendre
 (see Perron~\cite{oP13}, pp. 252-253):
\vspace{5pt}
\newline
\textbf{Tietze's Criterion}:

Let $ \{\,a_{n}\}_{n = 1}^{\infty}$ be a sequence of integers and
 $\{\,b_{n}\}_{n = 1}^{\infty}$
be a sequence of positive integers, with
$\,a_{n} \, \neq \, 0$ \,for any $n$.
If there exists a positive integer  $N_{0}$  such that
\begin{equation}\label{E:Ttze1}
\begin{cases}
 b_{n} \geq |a_{n}|\, & \\
 b_{n} \geq |a_{n}|\,+ 1, & for\,\, a_{n+1}\,<\,0,
\end{cases}
\end{equation}
for all \,$n\, \geq N_{0}$  then \,
$K_{n = 1}^{\infty}\displaystyle{a_{n}/b_{n}}\,$
converges and it's limit is irrational.

Here we use the standard notations
\[
\K{n = 1}{N}\displaystyle{\frac{a_{n}}{b_{n}}} \,:=
\cfrac{a_{1}}{b_{1} + \cfrac{a_{2}}{b_{2} + \cfrac{a_{3}}{b_{3} +
\ldots +  \cfrac{a_{N}^{}}{b_{N} }}}}
  \]
\[
\phantom{asd} =  \frac{a_{1}}{b_{1}+}\, \frac{a_{2}}{b_{2}+}\,
\frac{a_{3}}{b_{3}+}\, \dots  \frac{a_{N}}{b_{N}}.
\]
We write $\displaystyle{A_{N}/B_{N}}$ for the above finite continued
fraction written as a rational function of the variables
$a_{1},...,a_{N},b_{1},...,b_{N}$. $A_{N}$ is the \emph{$N$-th canonical numerator}
and $B_{N}$ is the \emph{$N$-th canonical denominator}.

It is elementary that the $A_{N}$ and $B_{N}$ satisfy the following
recurrence relations:
\begin{align}\label{E:recur}
A_{N}&=b_{N}A_{N-1}+a_{N}A_{N-2},\\
B_{N}&=b_{N}B_{N-1}+a_{N}B_{N-2} \notag.
\end{align}

By $K_{n = 1}^{\infty}\displaystyle{a_{n}/b_{n}}$
 we  mean the limit
of the sequence $\{\displaystyle{A_{n}/B_{n}}\}$ as \,$n$\,goes to infinity,
if the limit exists.

This type of continued fraction
seems to have  been studied systematically for the first time in Perron~\cite{oP13},
 where  degrees through two for $a_n$
and degree one for $b_n$\, are studied. Lorentzen
and Waadeland~\cite{LW92} also study these cases in detail and they evaluate
all such continued fractions in terms of hypergeometric series.
There is presently no such systematic treatment for cases of
higher degree.

Ramanujan also gave several
 continued fraction identities in which the free parameters can be specialized
to integers to give the values of polynomial continued fractions. An example of such a continued fraction
is the following entry from Chapter 12 in the second notebook, (see \cite{B89}, page 119 ):

\vspace{10pt}

\textbf{Entry 13.}
\emph{
Let $a$, $b$ and $d$ be complex numbers such that either $d \not = 0$, $b \not = - k d$, where
$k$ is a nonnegative integer, and Re$((a-b)/d)<0$, \footnote{
This condition is incorrectly given as "Re$((a-b)/d)>0$" in  \cite{B89}.
}
or $d \not = 0$ and $a=b$, or $d=0$ and $|a|<|b|$. Then}
{\allowdisplaybreaks
\begin{equation*}
a=\frac{ab}{a+b+d}
\-
\frac{(a+d)(b+d)}{a+b+3d}
\-
\frac{(a+2d)(b+2d)}{a+b+5d}
\-
\cds.
\end{equation*}
}

\vspace{10pt}

 In \cite{BML02}, the first author and Douglas Bowman examined certain
classes of polynomial continued fractions in which the $a_{n}$ and $b_{n}$ were of arbitrarily high degree.
One phenomenon we observed was the existence of infinite families of in-equivalent polynomial
continued fractions in which each continued fraction converged to the same limit.
To illustrates this phenomenon, we cite the following example from  \cite{BML02}:
\begin{example}\label{C:c12}
Let $f_{n}$ be a polynomial in
$n$ such that $f_{n}\geq 1$ for  $ n\geq 1$ and let $m$ be a
positive integer. Then
\[
\K{n =1}{\infty}
\displaystyle{\frac{f_{n}((n^{2}+3n+2)nm+1)+ 2mn^{2}+6mn+4m-1 }
{f_{n}((n^{2}-1)nm+1) +2(n^{2}-1)m-2}}
= 6m+1.
\]
\end{example}

One reason to study  infinite families of polynomial continued fractions in which each continued fraction
converges to the same limit is that hopefully one might arrive at a classification of a "base
set"  such that every convergent polynomial
continued fraction has the same limit as one of
the continued fractions in the base set. One might further hope that every convergent polynomial
continued fraction could be shown to have the same limit as some polynomial
continued fraction which is at
most quadratic in the numerator and at most linear
 in the denominator, since the limits of  these types
of continued fraction can be determined as ratios of hypergeometric series.
It would be interesting
to see either a proof of this or to produce a counter-example.

In  \cite{BML02}, our results were derived from a theorem of
Pincherle (Theorem \ref{T:P*} below) and a variant of the Euler
transformation (see Equation \ref{E:eul}). In this present paper,
we study this phenomenon in greater depth, considering not only
further applications of Pincherle's theorem and the Euler formula,
but also an extension of Euler's formula, an extension of the
corresponding transformation that equates infinite products and
continued fractions, extensions and contractions of continued
fractions and  the Bauer-Muir transformation.

As special cases  of our results, we give the following examples
(proofs are found throughout the paper).

\vspace{5pt}

\textbf{Example} \ref{expinch3}.
Let $c_{n}$  be any polynomial in $n$ such that   $c_{n} \geq 2$
for $n\geq -1$. Then
\begin{multline*} 1=\frac{{c_0} +
{{c_1}}^2}{{{c_1}}^2} \+
  \frac{{c_0}\,\left( {c_1} + {{c_2}}^2 \right) }{{{c_2}}^2}
\+
\cds
\+
  \frac{{c_{n-2 }}\,\left( {c_{n-1 }} + {{c_n}}^2 \right) }{{{c_n}}^2}
\+
\cds .
\end{multline*}

\vspace{5pt}

\textbf{Example} \ref{ex:pi}.
Let  $A$ be a positive integer. Then
{\allowdisplaybreaks
\begin{multline*}
\frac{\pi}{4}
= \frac{1}{A}
 +
\frac{1}{1}
\-
\frac{4+A}{4-2A}
\+\\
\displaystyle{K_{n=3}^{\infty}}
\frac{(2n-3)(2n-5)
  [A(2n-7)+4(n-3)  ]
 \left [ A(2n-3)+4(n-1) \right ]
}
{2A(2n-5)+4(2n-3)}.
\end{multline*}
}

\vspace{5pt}

\textbf{Example} \ref{ex:zeta}.
If $\zeta(n)$ is the Riemann Zeta function and $A$ is a positive integer, then
{\allowdisplaybreaks
\begin{multline*}
\zeta (11) =
\frac{1}{A}
+
\frac{2A-1}{2A}
\-
\frac{2A(3A-2^{10})}{3A(1+2^{11})-2^{12}}\\
\+
\displaystyle{K_{n=3}^{\infty}}
\frac{-n(n-1)^{21}\left[A(n-1)-(n-2)^{10}\right]\left[A(n+1)-n^{10}\right]}
{A(n+1)\left[(n-1)^{11}+n^{11}\right]-2n^{11}(n-1)^{10}}.
\end{multline*}
}

\vspace{5pt}

\textbf{Example} \ref{binomex}.
If $A$ is any integer different from $-7$, then
\begin{multline*}
\sqrt [5] {\frac{12}{7}} =
\frac{A+7}{7}
\+
\frac{7(11A-7)}
{-4A+56}
\+\\
K_{n=3}^{\infty}
\frac{7(5n-11)(n-2)\left[(12n-37)A-7\right]\left[(12n-13)A-7\right]}
{-2A(12n^2-31n+16)+14(2+n)}.
\end{multline*}

\vspace{5pt}

\textbf{Example} \ref{sinprod}.
For each integer $A$,
\begin{multline*}
\frac{3 \sqrt{3}}{2 \pi} = 1+A+ \frac{-10A-2}{18} \+ \frac{2736
A+432}{-692 A-132} \+ \\
\underset{n=3}{\overset{\infty}{K}}
 \frac{
\begin{matrix}
9n(1-n)(3n-
2)(3n-4)\phantom{asddgfdgfdgsdfdsadadsas} \\
\phantom{asddgfdgfdgsdf} \times [(1+9A)n+A+1][(1+9A)n-17A-1]
\end{matrix}
} 
{-8\,\left( 1 + A \right) +
  2\,\left( 5 + 9\,A \right) \,n + 144\,A\,n^2 -
  18\,\left( 1 + 9\,A \right) \,n^3 }.
\end{multline*}

\vspace{5pt}

\textbf{Example} \ref{eex}.
If $A$ is a non-negative integer,
then
\begin{multline*}
e=2+
\frac{1}{1+A}
\+
\frac{1-2A(1+A)}{2(1+A)}
\+
\frac{2(1-3A(1+A))}{3-5A-6A^2}\\
\+
K_{n=4}^{\infty}
\frac{(n-1)\left [ 1-nA(1+A)\right ]\left [ 1-(n-2)A(1+A)\right ]}
{n-(n(n-1)-1)A-n(n-1)A^{2}}.
\end{multline*}

We begin with some applications of Pincherle's Theorem.

\section {Pincherle's Theorem}

\begin{theorem}\label{T:P*}
$($Pincherle \cite{P94}$)$  Let $ \{\,a_{n}\}_{n = 1}^{\infty}$,
$\{\,b_{n}\}_{n = 1}^{\infty}$ and
 $ \{\,G_{n}\}_{n = -1}^{\infty}$ be \\ sequences of real or complex
  numbers satisfying  $a_{n} \neq 0$ for $n\geq 1$
 and for all \,\, $n \geq 1$,
\begin{equation}\label{E:pn}
 G_{n} =   a_{n}G_{n-2}  + b_{n} G_{n - 1}.
\end{equation}

\vspace{5pt}

Let $\{B_{n}\}_{n=1}^{\infty}$
denote the denominator convergents of the continued fraction
 $\K{n = 1}{\infty}\displaystyle{\frac{a_{n}}{b_{n}}}$.
If $\lim_{n \to \infty}\displaystyle{G_{n}/B_{n}} = 0$
then  $\K{n = 1}{\infty}\displaystyle{\frac{a_{n}}{b_{n}}}$
converges and  its limit is $-G_{0}/G_{-1}$.
\end{theorem}

In \cite{BML02} we looked for solutions where $a_{n}$ and $b_{n}$   were polynomials. We
then showed that the corresponding $B_{n}$ and $G_{n}$ in the statement of the theorem
satisfied $\lim_{n \to \infty}\displaystyle{G_{n}/B_{n}} = 0$
 by using some easily deduced facts about the growth of $B_{n}$:

(i) \quad
Let $a_{n}$\, and \, $b_{n}$
be  non-constant polynomials in $n$ such that $a_{n} \geq 1$,
$b_{n} \geq 1$,
for $n\geq1$ and suppose $b_{n}$ is a polynomial of degree\,$k$.
If the leading coefficient of $b_{n}$ is $D$, then given $\epsilon >0$,
there exists
a positive constant $C_{1}=C_{1}(\epsilon) $ such that
$B_{n} \geq C_{1}(|D|/(1+\epsilon))^{n}(n!)^{k}$.

\vspace{5pt}

(ii) \quad
If $a_{n}$\, and \, $b_{n}$ are positive numbers $\geq 1$, then  there
exists a
positive
constant $C_{3}$ such that $B_{n} \geq C_{3}\phi^{n}$ for $n\geq 1$,
where $\phi$ is
the golden ratio\, $\displaystyle{(1 + \sqrt{5})/2}$.
\vspace{10pt}

The drawbacks to this approach are, firstly,  the difficulty if finding solutions to Equation \ref{E:pn}
 where
the $a_{n}$ and $b_{n}$ are polynomials in $n$ and, secondly, showing that
 $B_{n}$ and $G_{n}$ do indeed satisfy $\lim_{n \to \infty}\displaystyle{G_{n}/B_{n}} = 0$.

Here we describe a very simple way of using Pincherle's theorem to write down infinite
families of continued fractions such that each member of the family converges to the same limit.
Moreover, $\lim_{n \to \infty}\displaystyle{G_{n}/B_{n}}$ $ = $ $0$ is satisfied automatically and it is
not necessary to have $a_{n}$ and $b_{n}$ polynomials in $n$
 (Initially the $a_{n}$'s and $b_{n}$'s
are rational functions but $K_{n=1}^{\infty}a_{n}/b_{n}$ can be converted to a polynomial
continued fraction by a similarity transformation.

We first need some notation.  If $f_{n}=g_{n}/h_{n}$ is a rational function in $n$, where
$g_{n}$ and $h_{n}$ are polynomials in $n$, we define the \emph{degree} of $f_{n}$ to be the
degree of $g_{n}$ minus the degree of $h_{n}$ and we define the \emph{leading coefficient} of $f_{n}$
to be quotient of the leading coefficient of $g_{n}$ and the leading coefficient of $h_{n}$.

\begin{proposition}\label{P:Pinchprop}
Let $H_{n}$ and $b_{n}$ be rational functions in $n$ such that
$H_{n}>0$ for $n\geq-1$ and $b_{n}>0$ for $n \geq 1$. We further
assume that $b_{n}$ has degree greater than $0$ or, if it has
degree $0$, then its leading coefficient is greater than 1. Then
the continued fraction
{\allowdisplaybreaks
\begin{equation}\label{E:Pinch1}
K_{n=1}^{\infty} \frac{(H_{n}+b_{n}H_{n-1})/
                                         H_{n-2}}
                                 {b_{n}}
\end{equation}
}
 converges and its limit is $H_{0}/H_{-1}$.
\end{proposition}
\begin{proof}
Let $G_{n}=(-1)^{n}H_{n}$ and $a_{n}=(H_{n}+b_{n}H_{n-1})/ H_{n-2}$. The sequences
 $\{G_{n}\}_{n=-1}^{\infty}$, $\{a_{n}\}_{n=-1}^{\infty}$ and $\{b_{n}\}_{n=-1}^{\infty}$ are easily
seen to satisfy Equation \ref{E:pn}. Furthermore, the recurrence relations  at \eqref{E:recur}
and the conditions on the $b_{n}$ give that $B_{n}$ grows at least as fast as either  $D^{n}$, for some fixed
$D>1$ or  $(n!)^{\delta}$, for some fixed $\delta >0$. Since $G_{n}$ is of polynomial growth,
it follows that $\lim_{n \to \infty}\displaystyle{G_{n}/B_{n}} = 0$.
Thus the continued fraction at \eqref{E:Pinch1} converges and  its limit is $-G_{0}/G_{-1}=H_{0}/H_{-1}$.
\end{proof}
As an illustration, we have the following example.
\begin{example}\label{EX2}
Let $H_{n}=n+2$ and let $b_{n}$ be any polynomial in $n$ such that
$b_{n} \geq 2$ for $n \geq 1$. Then
{\allowdisplaybreaks
\begin{multline*}
2=
\frac{3+2 b_{1}}{b_{1}}
\+
\frac{4+3b_{2}}{2b_{2}}
\+
\frac{2(5+4b_{3})}{3b_{3}}
\+
\cds \\
\cds
\+
\frac{(n-1)(n+2+(n+1)b_{n})}{nb_{n}}
\+
\cds
.
\end{multline*}
}
\end{example}
This follows from Proposition \ref{P:Pinchprop}, after a similarity transformation
to clear denominators.

Remarks:
1) The restrictions that $H_{n}>0$ and $b_{n}>0$ are not so severe since, if we restrict to polynomial
continued fractions $K_{n=1}^{\infty}a_{n}/b_{n}$ in which the polynomials $a_{n}$ and $b_{n}$ have
leading positive coefficients, then some tail of the continued fraction will have all $a_{n}$, $b_{n}>0$
so that, if one can find the limit of the tail, the continued fraction $K_{n=1}^{\infty}a_{n}/b_{n}$ reduces
to a finite continued fraction.

2) The limit of the continued fraction is independent of $b_{n}$ so that one has an infinite
class of continued fractions converging to the same limit.

3) Although we restrict $H_{n}$ and $b_{n}$ to be rational functions of $n$,
the result can easily be seen
to be true for more general sequences of positive numbers,
$\{H_{n}\}$ and $\{b_{n}\}$, provided the $b_{n}$ satisfy
certain size or growth conditions.

4) Convergence to the stated limit may also hold for other polynomial sequences $\{b_{n}\}$
which do not satisfy the conditions of the proposition (For example, the continued fraction in
Example \ref{EX2} converges when $b_{n}$ is the constant polynomial $1$).

We next write the rational functions $H_{n}$ and $b_{n}$ in terms of the polynomials defining them
to obtain a result about polynomial continued fractions.  This allows greater flexibility in deriving
the limits of infinite families of polynomial continued fractions.

\begin{corollary}
Let $f_{n}$, $g_{n}$, $c_{n}$ and $d_{n}$ be polynomials in $n$
such that $f_{n}$, $g_{n} \not =0$ for $n\geq -1$, $f_{n}$ and
$g_{n}$   have the same sign for each  $n\geq -1$,  $c_{n}$,
$d_{n} \not =0$ for $n\geq 0$ and $c_{n}$ and  $d_{n}$   have the
same sign for each  $n\geq 0$. Suppose further that the degree of
$c_{n}$ is greater than the degree of $d_{n}$ or, if the degrees
are equal, that the leading coefficient of $c_{n}$ is greater than
the leading coefficient of $d_{n}$. Then the continued fraction
{\allowdisplaybreaks
\begin{multline}\label{E:PinchPoly}
\frac{{g_{-1}}\,\left( {d_1}\,{f_1}\,{g_0} + {c_1}\,{f_0}\,{g_1}
\right) }
   {{c_1}\,{f_{-1}}\,{g_0}\,{g_1}}
\+\frac{{d_1}\,{f_{-1}}\,{{g_0}}^2\,
     \left( {d_2}\,{f_2}\,{g_1} + {c_2}\,{f_1}\,{g_2} \right) }
{{c_2}\,{f_0}\,{g_2}}
\+
\\
  \frac{{d_2}\,{f_0}\,{g_1}\,\left( {d_3}\,{f_3}\,{g_2} + {c_3}\,{f_2}\,{g_3} \right) }
   {{c_3}\,{f_1}\,{g_3}}
\+
\frac{{d_3}\,{f_1}\,{g_2}\,
     \left( {d_4}\,{f_4}\,{g_3} + {c_4}\,{f_3}\,{g_4} \right) }
{{c_4}\,{f_2}\,{g_4}}
\\
\+
\cds
  \frac{{d_{n-1}}\,{f_{ n-3}}\,{g_{n-2}}\,
     \left( {d_n}\,{f_n}\,{g_{ n-1}} + {c_n}\,{f_{n-1}}\,{g_n} \right) }{{c_n}\,
     {f_{ n-2}}\,{g_n}}
\+
\cds
\end{multline}
}
 converges and its limit is $f_{0}\,g_{-1}/(g_{0}\,f_{-1})$.
\end{corollary}
\begin{proof}
In Proposition \ref{P:Pinchprop}, let $H_{n}=f_{n}/g_{n}$ and $b_{n}=c_{n}/d_{n}$.
The continued fraction at \ref{E:Pinch1} is equivalent to the continued fraction at
 \ref{E:PinchPoly}, after a similarity transformation. The conditions of Proposition \ref{P:Pinchprop}
are satisfied  and the result follows.
\end{proof}
Remark: The limit is independent of the polynomials $c_{n}$ and $d_{n}$, provided they
satisfy the conditions of the corollary.

If we let $f_{n}=n^2+1$ and $d_{n}=g_{n}=1$ we get the following
example.
{\allowdisplaybreaks
\begin{example}\label{Ex1} Let
$c_{n}$  be any polynomials in $n$ such that   $c_{n} \geq 2$ for
$n\geq 0$. Then
{\allowdisplaybreaks
\begin{multline*}
\frac{1}{2}=
\\
\frac{{c_1} + 2}{2\,{c_1}}
\+
\frac{2\,\left( 2\,{c_2} + 5\, \right) }
   {{c_2}}
\+
\frac{\left( 5\,{c_3} + 10 \right) }{2\,{c_3}}
\+
  \frac{2\,\left( 10\,{c_4} + 17 \right) }{5\,{c_4}}
\+
\frac{5\,\left( 17\,{c_5} + 26 \right) }{10\,{c_5}}
\\
\+
\cds
  \frac{\left( 1 + {\left( -3 + n \right) }^2 \right) \,
     \left( \left( 1 + {\left( -1 + n \right) }^2 \right) \,{c_n} +
       \left( 1 + n^2 \right)  \right) }
{\left( 1 + {\left( -2 + n \right) }^2 \right)
       \,{c_n}}
\+
\cds .
\end{multline*}
}
\end{example}
}
 If we let $f_{n}=d_{n}=1$ and $g_{n}=c_{n}$ and cancel a factor
of $c_{-1}/c_{0}$ from the continued fraction and its limit, we
get the following result.
{\allowdisplaybreaks
\begin{example}\label{expinch3} Let $c_{n}$
be any polynomial in $n$ such that   $c_{n} \geq 2$ for $n\geq
-1$. Then
{\allowdisplaybreaks
\begin{multline*} 1=\frac{{c_0} +
{{c_1}}^2}{{{c_1}}^2} \+
  \frac{{c_0}\,\left( {c_1} + {{c_2}}^2 \right) }{{{c_2}}^2}
\+
\cds
\+
  \frac{{c_{-2 + n}}\,\left( {c_{-1 + n}} + {{c_n}}^2 \right) }{{{c_n}}^2}
\+
\cds .
\end{multline*}
}
\end{example}
}
 We next study a generalization of Euler's transformation of a
series into a continued fraction.

\section{The Euler Transformation of an Infinite Series}

 In 1775, Daniel Bernoulli \cite{B75}
proved the following result (see, for example, \cite{K63}, pp.
11--12).
{\allowdisplaybreaks
\begin{proposition}\label{pber} Let
$\{K_{0},K_{1}, K_{2},\ldots\}$ be a sequence of complex numbers
such that $K_{i}\not = K_{i-1}$, for $i=1,2,\ldots$. Then
$\{K_{0},K_{1}, K_{2},\ldots\}$  is the sequence of approximants
of the continued fraction
{\allowdisplaybreaks
\begin{multline}\label{ber1}
K_{0}+\frac{K_{1}-K_{0}}{1} \+ \frac{K_{1}-K_{2}}{K_{2}-K_{0}} \+
\frac{(K_{1}-K_{0})(K_{2}-K_{3})}
            {K_{3}-K_{1}}
\+\\
\ldots
\+
\frac{(K_{n-2}-K_{n-3})(K_{n-1}-K_{n})}
            {K_{n}-K_{n-2}}
\+
\ldots
.
\end{multline}
}
\end{proposition}
}
 In particular, if $\{K_{n}\}$ is a convergent sequence, one gets
a convergent continued fraction.

If one lets $K_{n}=\sum_{k=0}^{n}a_{k}$, one gets Euler's
transformation of a series into a continued fraction \cite{E27}:
{\allowdisplaybreaks
\begin{equation}\label{E:eul} \sum_{k =
0}^{n}a_{k} =
 a_{0} + \frac{a_{1}}{1}
\+
 \frac{-a_{2}}{a_{1}+a_{2}}
\+
\frac{-a_{1}a_{3}}{a_{2}+a_{3}}
\+
\cds
\+
\frac{-a_{n-2}a_{n}}{a_{n-1}+a_{n}}
\+
\cds .
\end{equation}
}
 For example, applying Euler's transformation to the well-known
series for $\pi/4$,
{\allowdisplaybreaks
\begin{equation}\label{piser} \frac{\pi}{4} =
\frac{1}{1} - \frac{1}{3} + \frac{1}{5} + \cdots +
\frac{(-1)^{n-1}}{2n-1} + \cdots ,
\end{equation}
}
 gives Lord Brouncker's continued fraction \eqref{E:lb}
 (after inversion and some similarity transformations to clear fractions).

We now give a generalization of the Euler transformation.
{\allowdisplaybreaks
\begin{proposition}\label{eug} Let $\sum_{k =
0}^{\infty}a_{k}$ be a convergent series  and let $\{b_{n}\}$ be
any sequence whose limit is zero such that $a_{n}+b_{n}-b_{n-1}
\not = 0$, for $n\geq 1$. Then
{\allowdisplaybreaks
\begin{multline}\label{eu}
\sum_{i=0}^{\infty}a_{i} = a_{0}+
b_{0}+\frac{a_{1}+b_{1}-b_{0}}{1} \+
\frac{-a_{2}+b_{1}-b_{2}}{a_{2}+a_{1} +b_{2}-b_{0}}\\
\+
K_{n=3}^{\infty}\frac{-(a_{n-2}+b_{n-2}-b_{n-3})(a_{n}+b_{n}-b_{n-1})}
{a_{n}+a_{n-1}+b_{n}-b_{n-2}}.
\end{multline}
}
\end{proposition}
}
\begin{proof}
This follows immediately from Proposition \ref{pber}, upon letting $K_{n}=\sum_{k=0}^{n}a_{k}+b_{n}$,
noting that $\lim_{n \to \infty}K_{n} = \sum_{i=0}^{\infty}a_{i}$.
\end{proof}
For our first example, we consider the series for $\pi/4$ above \eqref{piser}.
\begin{example}\label{ex:pi}
Let  $f_{n}$ be any polynomial which is positive for $n\geq 1$ and, for ease of notation,
define $g_{n} := f_{n}f_{n-1}+(2n-1)(f_{n}+f_{n-1})$. Then
{\allowdisplaybreaks
\begin{multline*} \frac{\pi}{4} =
\frac{-1}{f_{0}} + \frac{g_{1}}{f_{1}f_{0}} \+
\frac{f_{0}^2g_{2}}{2f_{2}f_{0}+3(f_{2}-f_{0})}\\
\+
K_{n=3}^{\infty}
\frac{(2n-3)^2g_{n-2}g_{n}}
{2f_{n}f_{n-2}+(2n-1)(2n-3)(f_{n}-f_{n-2})}.
\end{multline*}
}
\end{example}
This follows from setting  $b_{n}=(-1)^{n-1}/f_{n}$ in Proposition
\ref{eug} and simplifying the continued fraction. If we let
$f_{n}=A(2n-1)$, where $A$ is a positive integer, we get, after
some further simplification, that
{\allowdisplaybreaks
\begin{multline*} \frac{\pi}{4} = \frac{1}{A}
 +
\frac{1}{1}
\-
\frac{4+A}{4-2A}
\+\\
K_{n=3}^{\infty}
\frac{(2n-3)(2n-5)
  [A(2n-7)+4(n-3)  ]
 \left [ A(2n-3)+4(n-1) \right ]
}
{2A(2n-5)+4(2n-3)}.
\end{multline*}
}

As a second example, we consider the series representation for
$\zeta (k)$, $k$ an integer greater than $1$:
{\allowdisplaybreaks
\begin{equation*}
\zeta(k)=\sum_{n=1}^{\infty}\frac{1}{n^{k}}.
\end{equation*}
} If the Euler transformation is applied to this series (letting
$a_{0}=0$ and $a_{i} =1/i^k$, for $i\geq 1$), one easily gets that
{\allowdisplaybreaks
\begin{equation*} \zeta(k)= \frac{1}{1} \-
\frac{1}{2^{k}+1} \- \frac{2^{2k}}{3^{k}+2^{k}} \-
\frac{3^{2k}}{4^{k}+3^{k}} \- \cds \-
\frac{(n-1)^{2k}}{n^{k}+(n-1)^{k}} \- \cds .
\end{equation*}
}
 If we set $b_{n}=1/d_{n}$ in Proposition \ref{eug}, where
$d_{n}$ is a polynomial of degree at least 1 such that
$g_{n}:=d_{n}d_{n-1}+n^{k}(d_{n-1}-d_{n}) \not = 0$ for $n \geq
1$, then
{\allowdisplaybreaks
\begin{multline*} \zeta (k) =
\frac{1}{d_{0}} + \frac{g_{1}}{d_{0}d_{1}} \-
\frac{d_{0}^2g_{2}}{d_{2}d_{0}(1+2^{k})+2^{k}(d_{0}-d_{2})}\\
\+
\displaystyle{K_{n=3}^{\infty}}
\frac{-(n-1)^{2k}g_{n-2}g_{n}}{d_{n}d_{n-2}((n-1)^k+n^k)+(n-1)^kn^k(d_{n-2}-d_{n})}.
\end{multline*}
}
 We can specialize further to get the following continued
fraction.
\begin{example}\label{ex:zeta}
If $A$ is a positive integer, then
{\allowdisplaybreaks
\begin{multline*} \zeta (k) = \frac{1}{A} +
\frac{2A-1}{2A} \-
\frac{2A(3A-2^{k-1})}{3A(1+2^k)-2^{k+1}}\\
\+
\displaystyle{K_{n=3}^{\infty}}
\frac{-n(n-1)^{2k-1}\left[A(n-1)-(n-2)^{k-1}\right]\left[A(n+1)-n^{k-1}\right]}
{A(n+1)\left[(n-1)^k+n^k\right]-2n^k(n-1)^{k-1}}.
\end{multline*}
}
\end{example}
This follows from defining $d_{n}:=A(n+1)$ and simplifying the continued fraction.

For our third example, we consider the binomial series. Let
$|x|<1$ and $\alpha \in \mathbb{R}$. Then
{\allowdisplaybreaks
\begin{align}\label{eqbinomser} (1+x)^{\alpha}
&=
1+ \alpha\, x + \frac{\alpha (\alpha-1)}{2!}\,x^2 + \frac{\alpha (\alpha-1)(\alpha-2)}{3!}\,x^3 + \cdots\\
&=\sum_{n=0}^{\infty}\frac{(\alpha)_{n}}{n!}\,x^{n}, \notag
\end{align}
}
 where $(\alpha)_{n}=\alpha(\alpha-1)\cdots(\alpha -n+1)$ denotes
the \emph{falling factorial}. If we let $b_{n} =
r_{n}x^{n}(\alpha)_{n}/n!$, where $r_{n}$ is any rational function
such that
\[
g_{n}:=(\alpha - n + 1)\,x\,(1+r_{n})-n r_{n-1} \not = 0, \hspace{34pt} n \geq 1,
\] then
Proposition \ref{eug} gives that
{\allowdisplaybreaks
\begin{multline}\label{bingeneq}
(1+x)^{\alpha} = 1+r_{0} + \frac{g_{1}}{1} \- \frac{\alpha \, x\,
g_{2}}{\alpha \, x\, \left [(\alpha - 1)x (1+r_{2})+2 \right ] -
2r_{0}} \+
\\
K_{n=3}^{\infty}
\frac{-(n-1)\, x \, (\alpha -n+2)g_{n-2}g_{n}}
{
(\alpha -n+2)\,x\,\left[(\alpha -n+1)x(1+r_{n})+n\right]-n(n-1)r_{n-2}
}.
\end{multline}
}
 If we specialize by  letting $\alpha=1/5$, $x=5/7$ and
$r_{n}=A\,n-1$, ($A \not = -7)$ and then simplifying the continued
fraction,  we get that
\begin{example}\label{binomex}
If $A$ is any integer different from $-7$, then
{\allowdisplaybreaks
\begin{multline*} \sqrt [5] {\frac{12}{7}} =
\frac{A+7}{7} \+ \frac{7(11A-7)} {-4A+56}
\+\\
K_{n=3}^{\infty}
\frac{7(5n-11)(n-2)\left[(12n-37)A-7\right]\left[(12n-13)A-7\right]}
{-2A(12n^2-31n+16)+14(2+n)}.
\end{multline*}
}
\end{example}
Remark: It follows from Equation \ref{bingeneq} that  every real number of the form $(p/q)^{r/s}$,
where $p$, $q$, $r$ and $s$ are integers can be expanded in infinitely many ways
as a polynomial continued fraction.

We next consider a generalization of the formula transforming an
infinite product to a continued fraction.

\section{The Transformation of  Infinite Products to  Continued Fractions }
In Equation \ref{ber1}, if  $K_{i} \not =0$ for $i \geq 1$,
then the continued fraction on the left side
can be re-written as
\begin{multline}\label{ber2}
K_{0}+
\frac{K_{1}-K_{0}}{1}
\-
\frac{K_{1}/K_{0}(K_{2}/K_{1}-1)}{K_{2}/K_{0}-1}\\
\-
\frac{K_{2}/K_{1}(K_{1}/K_{0}-1)(K_{3}/K_{2}-1)}{K_{3}/K_{1}-1}
\- \\
\cds
\-
\frac{K_{n-1}/K_{n-2}(K_{n-2}/K_{n-3}-1)(K_{n}/K_{n-1}-1)}{K_{n}/K_{n-2}-1}
\-
\cds
.
\end{multline}

In particular, if $K_{0}=1$ and $K_{n}= \prod_{i=1}^{n}a_{i}$ for $n \geq 1$, where
$ \prod_{i=1}^{\infty}a_{i}$ is
a convergent infinite product with no $a_{i}=0$ or $1$ , one has that
\begin{multline}\label{berprod}
\prod_{i=1}^{\infty}a_{i}
=
1 + \frac{a_{1}-1}{1}
\-
\frac{a_{1}(a_{2}-1)}{a_{2}a_{1}-1}
\+
K_{n=3}^{\infty}\frac{-a_{n-1}(a_{n-2}-1)(a_{n}-1)}{a_{n}a_{n-1}-1}.
\end{multline}
This transformation is not so well known as Euler's transformation
of an infinite series to a continued fraction. As with the Euler
transformation, it is easy to generalize the transformation at
\eqref{berprod}.
\begin{proposition}
Let $\prod_{i=1}^{\infty}a_{i}$ be a convergent infinite product with $a_{i} \not = 0$ for
$i\geq 1$ and let $\{b_{n}\}_{n=0}^{\infty}$ be any sequence
whose limit is $1$  such that $a_{i}b_{i}-b_{i-1} \not = 0$, for $i \geq 1$. Then
\begin{multline}\label{berprodg}
\prod_{i=1}^{\infty}a_{i}
=
b_{0} + \frac{a_{1}b_{1}-b_{0}}{1}
\-
\frac{a_{1}(a_{2}b_{2}-b_{1})}
{a_{2}a_{1}b_{2}-b_{0}}\\
\+
K_{n=3}^{\infty}
\frac{-a_{n-1}(a_{n-2}b_{n-2}-b_{n-3})(a_{n}b_{n}-b_{n-1})}{a_{n}a_{n-1}b_{n}-b_{n-2}}.
\end{multline}
\end{proposition}
\begin{proof}
This follows from Proposition  \ref{pber} upon setting $K_{0}=b_{0}$ and, for
$n \geq 1$, setting  $K_{n}$ $=$
$b_{n}$ $ \prod_{i=1}^{n}a_{i}$
and then simplifying the continued fraction.
\end{proof}
As an example, we consider the following  infinite product identity.
\begin{equation*}
\frac{\sin \pi x}{\pi x}=\prod_{n=1}^{\infty}\left ( 1-\frac{x^2}{n^2}\right).
\end{equation*}
Set $x=1/m$, where $m$ is a positive integer, and $b_{n}=1+A/(n+1)$,
where $A$ is an integer.
For ease of notation,
let 
{\allowdisplaybreaks
\begin{align*}
g_{n}&=(m^2n+1)A+(n+1),\\
h_{n}&= (m^2n^2-1)(m^2(n-1)^2-1)(A+n+1)\\
&\phantom{aaaaaaaaaaaaaaaaaaaaaaaaaaaaaaaa}-m^4n^2(n^2-1)(A+n-1).
\end{align*}
}
Proposition \ref{berprodg}  then gives, after simplifying the continued fraction, that
{\allowdisplaybreaks
\begin{multline*}
\frac{m \sin \pi /m}{\pi }=
1+A+
\frac{(m^2-1)(A+2)-2 m^2(A+1)}{2m^2}\\
\+
\frac{2m^2(m^2-1)g_{2}}
{h_{2}}
\+K_{n=3}^{\infty}
\frac{-n(n-1)m^2(m^2(n-1)^2-1)g_{n-2}g_{n}
}
{h_{n}}
.
\end{multline*}
} If we specialize further and let $m=3$, we have that
\begin{example}\label{sinprod}
For each integer $A$,
{\allowdisplaybreaks
\begin{multline*}
\frac{3 \sqrt{3}}{2 \pi}
=
1+A+
\frac{-10A-2}{18}
\+
\frac{2736 A+432}{-692 A-132}
\+\\
\underset{n=3}{\overset{\infty}{K}}
 \frac{
\begin{matrix}
9n(1-n)(3n-2)(3n-4) \phantom{sdadasdasdasdasddsadaassdsa} \\
 \phantom{sdadasdasdasda} \times \left[ (1+9A)n+A+1\right
] \left[ (1+9A)n-17A-1\right]
\end{matrix}
}
{(9n^2-1)(9(n-1)^2-1)(A+n+1)-81n^2(n^2-1)(A+n-1)}.
\end{multline*}
}
\end{example}

\section{Extensions, Contractions and the Bauer-Muir Transformation}
We start with the concepts of extensions and
contractions of continued fractions. Before coming to details, we borrow some notation from \cite{LW92}
(page 83).
A continued fraction $d_{0}+K_{n=1}^{\infty}c_{n}/d_{n}$ is said to be a \emph{contraction} of the
continued fraction $b_{0}+K_{n=1}^{\infty}a_{n}/b_{n}$ if its classical approximants
$\{g_{n}\}$ form a subsequence of the classical approximants $\{f_{n}\}$
 of $b_{0}+K_{n=1}^{\infty}a_{n}/b_{n}$. In this case $b_{0}+K_{n=1}^{\infty}a_{n}/b_{n}$
is called an \emph{extension} of $d_{0}+K_{n=1}^{\infty}c_{n}/d_{n}$.

We call $d_{0}+K_{n=1}^{\infty}c_{n}/d_{n}$  a \emph{canonical contraction} of
 $b_{0}+K_{n=1}^{\infty}a_{n}/b_{n}$ if
\begin{align*}
&C_{k}=A_{n_{k}},& &D_{k}=B_{n_{k}}& &\text{ for } k=0,1,2,3,\ldots \, ,\phantom{asdasd}&
\end{align*}
where $C_{n}$, $D_{n}$, $A_{n}$ and $B_{n}$ are canonical numerators and denominators
of $d_{0}+K_{n=1}^{\infty}c_{n}/d_{n}$ and $b_{0}+K_{n=1}^{\infty}a_{n}/b_{n}$ respectively.

From \cite{LW92} (page 83 and page 85)we also have the following theorems.
\begin{theorem}\label{T:t1}
The canonical contraction of $b_{0}+K_{n=1}^{\infty}a_{n}/b_{n}$ with
\begin{align*}
&C_{k}=A_{2k}& &D_{k}=B_{2k}& &\text{ for } k=0,1,2,3,\ldots \, ,&
\end{align*}
exists if and only if $b_{2k} \not = 0 for K=1,2,3,\ldots$, and in this case is given by
\begin{equation}\label{E:evcf}
b_{0}
+
\frac{b_{2}a_{1}}{b_{2}b_{1}+a_{2}}
\-
\frac{a_{2}a_{3}b_{4}/b_{2}}{a_{4}+b_{3}b_{4}+a_{3}b_{4}/b_{2}}
\-
\frac{a_{4}a_{5}b_{6}/b_{4}}{a_{6}+b_{5}b_{6}+a_{5}b_{6}/b_{4}}
\+
\cds .
\end{equation}
\end{theorem}
The continued fraction \eqref{E:evcf} is called the \emph{even} part of $b_{0}+K_{n=1}^{\infty}a_{n}/b_{n}$.

\begin{theorem}\label{odcf}
The canonical contraction of $b_{0}+K_{n=1}^{\infty}a_{n}/b_{n}$ with
$C_{0}=A_{1}/B_{1}$
\begin{align*}
&C_{k}=A_{2k+1}& &D_{k}=B_{2k+1}& &\text{ for } k=1,2,3,\ldots \, ,&
\end{align*}
exists if and only if $b_{2k+1} \not = 0 for K=0,1,2,3,\ldots$, and in this case is given by
\begin{multline}\label{E:odcf}
\frac{b_{0}b_{1}+a_{1}}{b_{1}}
-
\frac{a_{1}a_{2}b_{3}/b_{1}}{b_{1}(a_{3}+b_{2}b_{3})+a_{2}b_{3}}
\-
\frac{a_{3}a_{4}b_{5}b_{1}/b_{3}}{a_{5}+b_{4}b_{5}+a_{4}b_{5}/b_{3}}\\
\-
\frac{a_{5}a_{6}b_{7}/b_{5}}{a_{7}+b_{6}b_{7}+a_{6}b_{7}/b_{5}}
\-
\frac{a_{7}a_{8}b_{9}/b_{7}}{a_{9}+b_{8}b_{9}+a_{8}b_{9}/b_{7}}
\+
\cds .
\end{multline}
\end{theorem}
The continued fraction \eqref{E:odcf} is called the \emph{odd} part of $b_{0}+K_{n=1}^{\infty}a_{n}/b_{n}$.

One might expect that if a continued fraction is constructed whose
even part is the convergent continued fraction
$K_{n=1}^{\infty}a_{n}/b_{n}$, that one might have some
flexibility in the choice of the partial quotients of the extended
continued fraction. This is indeed the case and if, in addition,
the extended continued fraction can be constructed so that it too
converges, an infinite family of continued fractions with the same
limit as $K_{n=1}^{\infty}a_{n}/b_{n}$ can be found.

However, before discussing this we first consider another way of transforming a
continued fraction so as to produce an infinite family of continued fractions
with the same limit.

\vspace{5pt}

\textbf{Definition} (\cite{LW92}, page 76) The Bauer-Muir transform of a
continued fraction $b_{0}+K(a_{n}/b_{n})$ with respect to a sequence $\{w_{n}\}$
from $\mathbb{C}$ is the continued fraction $d_{0}+K(c_{n}/d_{n})$ whose canonical
numerators $C_{n}$ and denominators $D_{n}$ are given by
\begin{align}
&C_{-1}=1,& &D_{-1}=0,&\\
&C_{n}=A_{n}+w_{n}A_{n-1},& &D_{n}=B_{n}+w_{n}B_{n-1}& \notag
\end{align}
for $n=0,1,2 \ldots$, where $\{A_{n}\}$ and $\{B_{n}\}$ are the canonical numerators
and  denominators of $b_{0}+K(a_{n}/b_{n})$.

This transformation dates back to the 1870's and the work of Bauer
\cite{B72} and Muir \cite{M77}. From \cite{LW92}, page 76, there
is the following theorem:
\begin{theorem}
The Bauer-Muir transform of $b_{0}+K(a_{n}/b_{n})$ with respect to $\{w_{n}\}$
from $\mathbb{C}$ exists if and only if
\begin{equation}
a_{n}-w_{n-1}(b_{n}+w_{n}) \not = 0 \text{  for }n=1,2,3,\ldots .
\end{equation}
If it exists, then it is given by
\begin{equation}
b_{0}+w_{0}
+
\frac{a_{1}-w_{0}(b_{1}+w_{1})}{b_{1}+w_{1}}
\+
\frac{c_{2}}{d_{2}}
\+
\frac{c_{3}}{d_{3}}
\+
\cds
\end{equation}
where
\begin{align}
c_{n}&=a_{n-1}
\frac{a_{n}-w_{n-1}(b_{n}+w_{n})}{a_{n-1}-w_{n-2}(b_{n-1}+w_{n-1})}\,,\\
d_{n}&=b_{n}+w_{n}-w_{n-2}
\frac{a_{n}-w_{n-1}(b_{n}+w_{n})}{a_{n-1}-w_{n-2}(b_{n-1}+w_{n-1})}\,. \notag
\end{align}
\end{theorem}
We have changed the notation found in \cite{LW92} slightly.
One might expect that, by choosing
the sequence $\{w_{n}\}$ appropriately, one could construct an infinite family
of continued fractions with the same limit as the convergent continued fraction
$K(a_{n}/b_{n})$.

Curiously, these two methods
(extensions and contractions and the Bauer - Muir transform)
 of producing  infinite families
of continued fractions in which each continued fraction
has the same limit as the convergent
continued fraction $K(a_{n}/b_{n})$ are related. We believe the following
observation to be new.
\begin{theorem}
Let $\{w_{n}\}_{n=0}^{\infty}$ be a sequence from $\mathbb{C}$ such that
$w_{0}=0$ and $w_{n} \not = 0$ for $n \geq 1$.
We suppose further that $\{a_{n}\}_{n=1}^{\infty}$ and $\{b_{n}\}_{n=1}^{\infty}$
are sequences from $\mathbb{C}$  such that
$a_{n}-w_{n-1}(b_{n}+w_{n}) \not = 0 \text{  for }n=1,2,3,\ldots $.
Then the even
part of the continued fraction
{\allowdisplaybreaks
\begin{multline} \label{bmoe}
\frac{a_{1}}{b_{1}+w_{1} }
\+
\frac{-w_{1}}{1}
\+
\frac{a_{2}/w_{1}}{b_{2}+w_{2}-a_{2}/w_{1}}
\+
\frac{-w_{2}}{1}
\+
\frac{a_{3}/w_{2}}{b_{3}+w_{3}-a_{3}/w_{2}}
\+\\
\cds
\+
\frac{-w_{n-1}}{1}
\+
\frac{ a_{n}/w_{n-1} }{b_{n}+w_{n}-a_{n}/w_{n-1} }
\+
\frac{-w_{n}}{1}
\+
\cds
\end{multline}
}
is $K_{n=1}^{\infty}a_{n}/b_{n}$ and its odd part is equal to the
Bauer-Muir transform of $K_{n=1}^{\infty}a_{n}/b_{n}$ with
respect to the sequence $\{w_{n}\}_{n=0}^{\infty}$.
\end{theorem}
Remark: The theorem says that the odd part of
the continued fraction at \ref{bmoe} " is equal to"  the
Bauer-Muir transform of $K_{n=1}^{\infty}a_{n}/b_{n}$ with
respect to the sequence $\{w_{n}\}_{n=0}^{\infty}$ rather than
saying that it "is" this Bauer-Muir transform, since some transformations
need to be applied to the odd part to make it equal the stated
Bauer-Muir transform of $K_{n=1}^{\infty}a_{n}/b_{n}$.
\begin{proof}
That the even part of \eqref{bmoe} is $K_{n=1}^{\infty}a_{n}/b_{n}$
follows immediately from Theorem \ref{T:t1}. From Theorem \ref{odcf}, the odd
part of   \eqref{bmoe} is
{\allowdisplaybreaks
\begin{multline*}
\frac{a_{1}}{b_{1}+w_{1}}
-
\frac{a_{1}(-w_{1})
\displaystyle{
\frac{b_{2}+w_{2}-a_{2}/w_{1}}{b_{1}+w_{1}}}
}
{(b_{1}+w_{1})
\left [ \displaystyle{
\frac{a_{2}}{w_{1}}
}
+ \left (b_{2}+w_{2}-
\displaystyle{
\frac{a_{2}}{w_{1}}
}
\right ) \right ]
 -w_{1}(b_{2}+w_{2}-a_{2}/w_{1})}\\
\-
\frac{
a_{2}/w_{1}(-w_{2})
\displaystyle{
\frac{ b_{3}+w_{3}-a_{3}/w_{2}}{b_{2}+w_{2}-a_{2}/w_{1}}}
(b_{1}+w_{1})
}
{
a_{3}/w_{2}+(b_{3}+w_{3}-a_{3}/w_{2})-w_{2}
\displaystyle{
\frac{ b_{3}+w_{3}-a_{3}/w_{2}}{b_{2}+w_{2}-a_{2}/w_{1}}}
}\\
\-
\frac{
a_{3}/w_{2}(-w_{3})
\displaystyle{
\frac{ b_{4}+w_{4}-a_{4}/w_{3}}{b_{3}+w_{3}-a_{3}/w_{2}}
}
}
{
a_{4}/w_{3}+(b_{4}+w_{4}-a_{4}/w_{3})-w_{3}
\displaystyle{
\frac{
b_{4}+w_{4}-a_{4}/w_{3}}
{b_{3}+w_{3}-a_{3}/w_{2}}}
}
\-\\
\cds \\
\-
\frac{
a_{n-1}/w_{n-2}(-w_{n-1})
\displaystyle{
\frac{ b_{n}+w_{n}-a_{n}/w_{n-1}}{b_{n-1}+w_{n-1}-a_{n-1}/w_{n-2}}
}
}
{
\displaystyle{
\frac{a_{n}}{w_{n-1}}
}
+
\left (b_{n}+w_{n}-
\displaystyle{
\frac{a_{n}}{w_{n-1}}
}
\right )
-w_{n-1}
\displaystyle{
\frac{
b_{n}+w_{n}-a_{n}/w_{n-1}}
{b_{n-1}+w_{n-1}-a_{n-1}/w_{n-2}}}
}
\end{multline*}
}
$\phantom{asdasdasdadadsadsd}\-
\cds
$
{\allowdisplaybreaks
\begin{multline*}
= \frac{a_{1}}{b_{1}+w_{1}}
-
\frac{a_{1}
\displaystyle{
\frac{a_{2}-w_{1}(b_{2}+w_{2})}{(b_{1}+w_{1})^{2}}}}
{
\displaystyle{
\frac{
a_{2}+b_{1}(b_{2}+w_{2})}
{b_{1}+w_{1}}
}
}
\+
\frac{
a_{2}
\displaystyle{
\frac{a_{3}-w_{2}( b_{3}+w_{3})}{a_{2}-w_{1}(b_{2}+w_{2})}}
}
{
b_{3}+w_{3}-w_{1}
\displaystyle{
\frac{
a_{3}-w_{2}(b_{3}+w_{3})}
{a_{2}-w_{1}(b_{2}+w_{2})}}
} \\
\+
\frac{
a_{3}
\displaystyle{
\frac{a_{4}-w_{3}( b_{4}+w_{4})}{a_{3}-w_{2}(b_{3}+w_{3})}}
}
{
b_{4}+w_{4}-w_{2}
\displaystyle{
\frac{
a_{4}-w_{3}(b_{4}+w_{4})}
{a_{3}-w_{2}(b_{3}+w_{3})}}
}
\+\\
\cds
\+
\frac{
a_{n-1}
\displaystyle{
\frac{ a_{n}-w_{n-1}(b_{n}+w_{n})}{a_{n-1}-w_{n-2}(b_{n-1}+w_{n-1})}}
}
{
b_{n}+w_{n}-w_{n-2}
\displaystyle{
\frac{ a_{n}-w_{n-1}(b_{n}+w_{n})}
{a_{n-1}-w_{n-2}(b_{n-1}+w_{n-1})}}
}
\+
\cds .
\end{multline*}
}
Since $b_{0}=w_{0}=0$, all that is necessary to complete the proof  is to show that,
for $\alpha$ arbitrary,
{\allowdisplaybreaks
\begin{equation*}
\frac{a_{1}}{b_{1}+w_{1}}
\+
\frac{a_{2}-w_{1}(b_{2}+w_{2})}{b_{2}+w_{2}+\alpha}
=\frac{a_{1}}{b_{1}+w_{1}}
-
\frac{a_{1}
\displaystyle{
\frac{a_{2}-w_{1}(b_{2}+w_{2})}{(b_{1}+w_{1})^{2}}}}
{
\displaystyle{
\frac{
a_{2}+b_{1}(b_{2}+w_{2})}
{b_{1}+w_{1}} + \alpha}}.
\end{equation*}
}
This is immediate from the identity
{\allowdisplaybreaks
\begin{equation*}
\frac{E}{F}
\+
\frac{G}{H+\alpha}
=
\frac{E}{F}
-
\frac{EG/F^{2}}{(FH+G)/F+\alpha}.
\end{equation*}
}
\end{proof}
We next consider the question of when a continued fraction and its
Bauer-Muir transform converge to the same limit. One approach is to
apply general convergence theorems to the continued fraction at \eqref{bmoe}.
If this continued fraction converges, then its odd and even parts converge to the
same limit. Alternatively, we have the following proposition.
\begin{proposition}\label{propbmoe}
Let  $\{w_{n}\}_{n=0}^{\infty}$ be any rational function sequence
taking only positive values for $n\geq 1$ and suppose
$w_{0}=0$. Suppose further that the polynomial sequences
$\{a_{n}\}_{n=1}^{\infty}$ and $\{b_{n}\}_{n=1}^{\infty}$ take
only positive values for $n \geq 1$, that
$a_{n}-w_{n-1}(b_{n}+w_{n}) \not = 0 \text{  for }n=1,2,3,\ldots $
and that the continued fraction $K_{n=1}a_{n}/b_{n}$ converges.
Then the Bauer-Muir transform of $K_{n=1}a_{n}/b_{n}$ with respect
to $\{w_{n}\}_{n=0}^{\infty}$ converges to the same limit.
\end{proposition}
\begin{proof}
Suppose $K_{n=1}^{\infty}a_{n}/b_{n} = L$. Then
\[
\lim_{n \to \infty} L- \frac{A_{n}}{B_{n}}=\lim_{n \to \infty} \frac{A_{n+1}}{B_{n+1}}- \frac{A_{n}}{B_{n}}=0.
\]
and
{\allowdisplaybreaks
\begin{align*}
\left | L-\frac{A_{n+1}+w_{n}A_{n}}{B_{n+1}+w_{n}B_{n}} \right |
&\leq \left | L-\frac{A_{n}}{B_{n}} \right |
+
\left | \frac{A_{n+1}+w_{n}A_{n}}{B_{n+1}+w_{n}B_{n}}
-\frac{A_{n}}{B_{n}}
\right | \\
&= \left | L-\frac{A_{n}}{B_{n}} \right |
+
\left | \frac{A_{n+1}}{B_{n+1}}
-\frac{A_{n}}{B_{n}}
\right |
\frac{1}{1+w_{n}B_{n}/B_{n+1}}
\\
&\leq \left | L-\frac{A_{n}}{B_{n}} \right |
+
\left | \frac{A_{n+1}}{B_{n+1}}
-\frac{A_{n}}{B_{n}}
\right |.
\end{align*}
}
The last inequality follows since $w_{n}B_{n}/B_{n+1}\geq 0$.
Finally, let $n \to \infty $ to get the result.
\end{proof}
This can be restated as follows: If $a_{n}$ and $b_{n}$ are polynomials taking
only positive values for $n \geq 1$ and $w_{n}$ is a rational function
taking only positive values for $n \geq 1$ and $w_{0}$ is defined to be $0$, then
{\allowdisplaybreaks
\begin{align*}
K_{n=1}^{\infty}\frac{a_{n}}{b_{n}} =
\frac{a_{1}}{b_{1}+w_{1}}
\+
\frac{a_{1}(a_{2}-w_{1}(b_{2}+w_{2}))}{a_{1}(b_{2}+w_{2})}\\
\+K_{n=3}^{\infty}
\frac{
a_{n-1}\left[a_{n}-w_{n-1}(b_{n}+w_{n})\right ]
\left[a_{n-2}-w_{n-3}(b_{n-2}+w_{n-2})\right ]
}
{
a_{n-1}(b_{n}+w_{n})-w_{n-2}(a_{n}+b_{n-1}(b_{n}+w_{n}))
}. \notag
\end{align*}
}
As an example, we consider the well known continued fraction expansion for $e$:
{\allowdisplaybreaks
\begin{equation*}
e=2
+
\frac{2}{2}
\+
\frac{3}{3}
\+
\frac{4}{4}
\+
\cds
\+
\frac{n}{n}
\+
\cds
.
\end{equation*}
}
If we let $A$ be a non-negative integer and
define $w_{0}:=0$  and $w_{n}:=A(n+1)$ for $n \geq 1$, we get the example below.
\begin{example}\label{eex}
If $A$ is a non-negative integer,
then
{\allowdisplaybreaks
\begin{multline*}
e=2+
\frac{1}{1+A}
\+
\frac{1-2A(1+A)}{2(1+A)}
\+
\frac{2(1-3A(1+A))}{3-5A-6A^2}\\
\+
K_{n=4}^{\infty}
\frac{(n-1)\left [ 1-nA(1+A)\right ]\left [ 1-(n-2)A(1+A)\right ]}
{n-(n(n-1)-1)A-n(n-1)A^{2}}.
\end{multline*}
}
\end{example}

\section{Conclusion}
Let $\mathbb{P}$ denote the set of all polynomial
continued fractions.\footnote{Of course not every element of $\mathbb{P}$
converges to a real number.}  It is not difficult
to construct an injective map from $\mathbb{P}$ to $\mathbb{N}$ so that
$\mathbb{P}$ is a countable set and thus that almost all real numbers
do \emph{not} have a polynomial continued fraction expansion.
Let $\mathbb{P}'$ denote the set of real numbers which are the limits of
convergent continued fractions in $\mathbb{P}$.
Trivially, $\mathbb{Q} \in \mathbb{P}'$ and, as is shown by  examples from the
literature, many algebraic- and transcendental numbers also have a polynomial
continued fraction expansion.
We conclude with a number of questions about polynomial continued fractions which
we consider interesting.

\vspace{5pt}

1) Is there an equivalent classification of the set $\mathbb{P}'$, perhaps in terms
of the partial quotients in the regular continued fraction expansion of its elements
or in terms of irrationality measure?

\vspace{5pt}

2) Does every real algebraic number belong to $\mathbb{P}'$? If not, exhibit
a counter-example. Note that the answer to this second question is yes if 
every real algebraic number $\alpha$ can be expanded 
as a hypergeometric series $\alpha = \sum_{i=0}^{\infty}b_{i}$, where 
each $b_{i} \in \mathbb{Q}$ and $b_{n+1}/b_{n}= r(n)$, 
where $r(x) \in \mathbb{Q}(x)$. We are not aware 
if this question has been answered in the literature.

\vspace{5pt}

3) Since almost all real numbers do \emph{not} belong to $\mathbb{P}'$, exhibit
an example of such a number, or perhaps an infinite family of such numbers.

\vspace{5pt}

4) With a little more effort, Proposition \ref{propbmoe} can be extended to show that,
given any convergent polynomial continued fraction, there is an infinite family of
convergent polynomial continued fractions with the same limit. Let $\mathbb{P}_{n,d}$
denote the set of all polynomial continued fractions
 of degree at most $n$ in the numerator and degree at most $d$ in the denominator.
Does there exist a pair of non-negative integers $n$ and $d$ such that
every element of $\mathbb{P}$ that converges has the same limit as some element of
$\mathbb{P}_{n,d}$? Is this true for $n=2$ and $d=1$?
\footnote{This would be of interest since every element of $\mathbb{P}_{2,1}$
can be evaluated in terms of hypergeometric functions.}

\vspace{5pt}

5)
Trivially, if $\alpha \not = 0$ and $\alpha \in \mathbb{P}'$, then
$1/\alpha \in \mathbb{P}'$. Does $\mathbb{P}'$ have any further algebraic structure?
Is it true that if $\alpha$, $\beta \in \mathbb{P}'$, then $\alpha + \beta \in \mathbb{P}'$?
 Is it true that if $\alpha$, $\beta \in \mathbb{P}'$,
then $\alpha \times \beta \in \mathbb{P}'$? Even negative answers to these questions
would be of interest but it would probably not be easy to find  counter-examples.

\vspace{5pt}

\allowdisplaybreaks{}

\end{document}